\documentclass[leqno,12pt]{amsart} 
\usepackage{amssymb}
\usepackage{amsmath}
\usepackage{enumerate}
\usepackage{amsfonts}
\usepackage{color}
\usepackage{graphicx}
\usepackage{setspace}
\usepackage{verbatim}
\usepackage{url}
\newtheorem{thm}{Theorem}[section]

\newtheorem{prop}[thm]{Proposition}
\newtheorem{cor}[thm]{Corollary}

\theoremstyle{remark}

\theoremstyle{remark}

\theoremstyle{definition}
\newtheorem{df}[thm]{Definition}
\newtheorem{ex}[thm]{Example}

\numberwithin{equation}{section}

\newcommand{\G}{\mathcal{G}}

\newcommand{\R}{\mathbb{R}}

\newcommand{\T}{\mathbb{T}}

\def\8{\infty}
\def\blab#1{\begin{equation}\label{#1}}
\def\elab{\end{equation}}

\def\blab#1{\begin{equation}\label{#1}}
\def\elab{\end{equation}}

\begin{document}
\title[PC matrix decomposition]
{Pairwise Comparisons Matrix Decomposition into Approximation and Orthogonal Component Using Lie Theory} 

\author[Koczkodaj]{W.W. Koczkodaj}
\address{Computer Science\\ Laurentian University\\ Sudbury, Ontario P3E 2C6, Canada} 
\email{ wkoczkodaj@cs.laurentian.ca}

\author[Marek]{V.W. Marek}
\address{Department of Computer Science\\
	University of Kentucky\\ Lexington, KY 40506, USA}
\email{marek@cs.uky.edu}

\author[Yayli]{Y. Yayli}
\address{Department of Mathematics, Faculty of Science, Ankara University\\
TR-06100 Ankara, Turkey}
\email{Yusuf.Yayli@science.ankara.edu.tr} 

\begin{abstract}
		
This paper examines the use of Lie group and Lie Algebra theory to construct the geometry of pairwise comparisons matrices.  The Hadamard product (also known as coordinatewise, coordinate-wise, elementwise, or element-wise product) is analyzed in the context of inconsistency and inaccuracy by the decomposition method. 

The two designed components are the approximation and orthogonal components. The decomposition constitutes the theoretical foundation for the multiplicative pairwise comparisons. \\
		
\end{abstract}
\maketitle
	
\section{Introduction}

Pairwise comparisons (PCs) take place when we somehow compare two entities (objects or abstract concepts).  According to \cite{HHPRD}, Raymond Llull is credited for the first documented use of pairwise comparisons in ``A system for the election of persons'' (\textit{Artifitium electionis personarum}) before 1283 and in ``An electoral system'' (\textit{De arte eleccionis}) in 1299. Both manuscripts were handwritten (there was no scientific publication process established yet) for deciding the winner of elections.
	
There are two variants of pairwise comparisons: multiplicative and additive. The multiplicative PCs variant reflects a relationship:
$$A~is~x~times~[comparison~operator]~than~B$$ 
The additive type expresses: \textit{``by how much (the percentage is often used) one entity is [comparison operator]
than another entity''}. The comparison operator could be: \textit{bigger, better, more important}, or similar comparisons (see \cite{LDX2019}).

The multiplicative pairwise comparison is determined by the ratio of two entities. For instance, the constant $\pi$, is one of the most recognized ratios in mathematics is defined as the ratio of a circle circumference to its diameter. 

Entities may be physical objects or abstract concepts (e.g., software dependability or software safety). In both cases, we can provide an estimate based on expert assessment. However, physical measurements (e.g., area or weight) should be considered if they are available.

In practice, multiplicative (or ratio) PCs are more popular than additive PCs. However, they are more mathematically challenging than additive comparisons. The additive PCs can be produced from the multiplicative form by a logarithmic mapping (introduced in \cite{RAND}). Values of additive pairwise comparisons could be both negative and positive real numbers.
The logarithmic mapping is used to provide the proof of
inconsistency convergence and for the interpretation of the limit of convergence in \cite{HK1996, KS2010, KS2016}. Fuzzy extensions of pairwise comparisons are presented in \cite{D1999, KLH2006}. 

Pairwise comparisons are usually represented by a \textit{PC matrix}. In the case of multiplicative PCs, it is a matrix of ratios of entities with 1s on the main diagonal (for the entity being compared to itself) and reciprocal ($x$ and $1/x$) values in upper/lower triangles  as it is also reasonable to assume that
the ratio of $B/A$ is the reciprocal value of $A/B$. By the definition, ratios are strictly positive real numbers.

For entities which are abstract concepts (e.g., software quality and software safety), the division operation is undefined but using the ratio still makes sense (e.g., by stating: ``three times more important''). For this reason, ratios are often given by expert assessments.

The main goal of PCs is to split 1 into $n$ real values assigned to $n$ entities $E_i$, $i=1, \ldots, n$. We call them {\em weights}.


\subsection{Problem outline}
\label{outline}
		

In \cite{KSS2020} (a follow-up to \cite{KO1997}), Smarzewski has observed that PC matrices form a group under Hadamard product. In this paper, we show that it is, in fact, a Lie group. It is the first application of Lie theory to pairwise comparisons method.

\subsection{Group theory applications in Computer Science}

In pairwise comparisons, the use of abelian groups took place in \cite{CavDap09}. More recently, abelian groups were used in \cite{Ram16,CB2018,KMRS2019}.

The intensive Internet search has identified the first  publication related to Lie group as \cite{BK1972}.
Most computer science publications, which use Lie theory are related to graphics and vision (see \cite{BT2019,EI2020,KCD2009,KF2020,ZZH2019}), Machine Learning (see \cite{DAO2018,KF2020,LG2006}), and network (see \cite{ZCX2019}).



%


%



---

\subsection{Structure of this paper}	
A brief introduction to the pairwise comparisons method is provided in Section~\ref{case}. Section~\ref{case} also includes a simple example of using the PC method for
a generic exam.  Section~\ref{basic-def} provides reasoning for the necessity of theorems and propositions. Section~\ref{LieGroupIntro} shows the
construction of Lie groups and Lie algebra for PC matrices.
Section~\ref{exp-trans} introduces the exponential transformation and its properties. Section~\ref{dir-prod} presents the main theorem. Section~\ref{genSec} outlines the generalization of our results.

\section{Example of using pairwise comparisons in practice}
\label{case}
		
A Monte Carlo experiment for pairwise comparisons accuracy improvement was presented in \cite{K1996,K1998}. It provided statistical evidence that the accuracy gain was substantial by using pairwise comparisons. For simplicity, let us assume that
we have three entities to compare: $ A, B,  $ and $ C $.
The three comparisons are: $ A $ to $B$, $ A $ to $ C $, and $ B $ to $ C $. We assume the reciprocity of PC matrix $M$: $m_{ji}=1/m_{ij}$ which is reasonable (when comparing $ B $ to $ A $, we expect to get the inverse of  the comparison $ A $ to $ B $).  The exam is hence represented by the following PC matrix $M$:
	
\begin{equation} 
	M=[m_{ij}]=
	\begin{bmatrix}
	  1   & A/B &  A/C \\
	  B/A &  1  &  B/C \\
	  C/A & C/B &   1  \\
	\end{bmatrix}
\end{equation}
	
\noindent $A/B$ reads ``the ratio between A and B'' and may not necessarily be a result of the division (in the case of the exam problem, the use of division operation makes no mathematical sense but using the ratio is still valid).
	
\begin{figure}[h]
		\centering
		\includegraphics[width=1\linewidth]{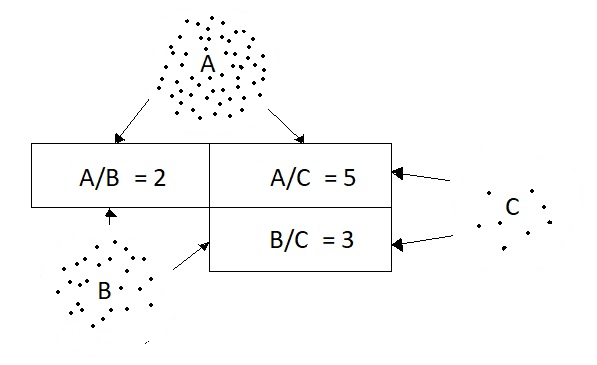}
		\caption[cycle]{An inconsistency indicator cycle}
		\label{fig:iicycle}
\end{figure}

Ratios of three entities create a {\em triad}  $(A/B, A/C,B/C)$. This triad is said to be \textit{consistent} provided $A/B*B/C=A/C$. It is illustrated by
Fig.~\ref{fig:iicycle}. Random numbers of dots are hard to count but can be compared in pairs as two random clouds. $[A/B]$ reflects the \textit{assessed} ratio of dots by expert opinions. 
	
A large enough number of dots represents the concept of (\textit{numerosity}). They may, for example, represent votes of experts. In an emergency situation (e.g., mine rescue), it is impossible to count votes in a short period of time. The exact number of votes is there but all we need is to assess the numerosity of votes .
	
Symbolically, in a PC matrix $M$, each triad (or a cycle) is defined by $(m_{ik},m_{ij},m_{kj})$. Such triad is {\em consistent} providing $(m_{ik}*m_{kj}=m_{ij})$. When all triads are consistent (known as the \textit{consistency condition} or \textit{transitivity condition}), the entire
PC matrix is considered \textit{consistent}.
	
	Looking at the above exam grading case, we have discovered a pairwise
comparisons method which can be used to construct a PC matrix. The
\textit{solution} to the PC matrix is a vector of weights which are geometric means of rows(for more detail see, for example, \cite{KK2015}). We usually normalize it to 1 as the sum. The justification for the use of the vector $w=[v_i]$ of geometric means (GM) of rows is not trivial and it is the subject of this paper. The exact reconstruction of the PC matrix (say $M$) via $M=[v_i/v_j]$ is guaranteed only for the consistent matrices.
	
In our example, the computed weights (as normalized geometric means of rows) are approximated to: $[0.58, 0.31, 0.11]$. By looking at the result, we can conclude that problem $A$ is the most difficult with the weight 0.58. The easiest problem is $C$ with the weight $0.11$.

One of the challenges of pairwise comparisons is the inconsistency of assessments. It is well demonstrated by Fig.~\ref{fig:iicycle}.  It seems that a trivial mistake took place: 6 should be in place of 5 since 2*3 gives
this value. However, there are no grounds to assume that 2 and 3 are accurate assessments. No specific accuracy assumption is made of any assessment.
		
\section{The problem formulation}
\label{basic-def}

This paper examines the use of group theory to construct the geometry of pairwise comparisons matrices and improve the consistency-driven process in \cite{JK1996, BFK2011}. The Hadamard product (also known as coordinatewise, coordinate-wise, elementwise, or element-wise product) is examined in the context of inconsistency and inaccuracy. To achieve this goal, we provide a proof that PC matrices are represented by a Lie group. Subsequently, a Lie algebra of
the Lie group of PC matrices is constructed. A decomposition method of PC matrices is introduced for the Hadamard product. One of the components is an approximation PC matrix and the other orthogonal component is interpreted as the approximation inaccuracy. The importance of selecting PC matrix components is also provided in this paper. Subgroups of the PC matrix Lie group have been identified and presented as an internal direct  product.

\section{Lie groups and Lie algebras of PC matrices}
\label{LieGroupIntro}

The monograph \cite{Tu2010} stipulates that, ``Intuitively, a manifold is a generalization of curves and surfaces to higher dimensions. It is locally Euclidean in that every
point has a neighborhood, called a chart, homeomorphic to an open subset of $\R^n$''. We find the above justification to be sufficient to be followed by computer science researchers.

A group that is also a differentiable (or smooth) manifold is called Lie group (after its proponent Sophus Lie).  According to \cite{AA2008}, a Lie group is an abstract group $\G$ with a smooth structure, that is: 
\begin{df}  \
\begin{enumerate} 
\item  $ \G $ is a group,
\item  $ \G $ is a smooth manifold,
\item  the operation $\G \times \G \longrightarrow \G,$ $(x,y) \longrightarrow
xy^{-1}$ is smooth.  
\end{enumerate}
\end{df}

\noindent Matrix Lie group operates on matrices.

\begin{df}
	The Lie algebra of a Lie group $ \G $ is the vector space $ \T_e\G $ equipped with the Lie bracket operation $ [~ ,~ ]$ of vector fields.
\end{df}

\noindent The bracket operation $ [~ , ~]$ is assumed to be bilinear, antisymmetric, and satisfies the Jacobi identity:  Cyclic$([X, [Y, Z]]) = 0$ for all $X, Y, Z$
belonging to this algebra.

%

\noindent Lie group and Lie algebra have been analyzed in 
\cite{YY2010,YY2013, YY2014,YY2016}.

\begin{prop}
For every dimension $n > 0$, the following group:
$$ \G=\{\  M=[m_{ij}]_{n \times n} | M \cdot M^T = I, m_{ij} = \frac{1}{m_{ji}}
> 0 \text{ for every } i,j = 1, 2,\ldots,n\}\ $$ 
is an abelian group of $n \times
n$ PC matrices with an operation
\begin{equation*}
\begin{split}
& \cdot : \G \times \G \longmapsto \G \\
& \mathit{defined \ by}\ \quad (M,N) \longrightarrow M \cdot N =[m_{ij} \cdot n_{ij}] 
\end{split}
\end{equation*}
where "$\cdot$" is the Hadamard product.
\end{prop}
\begin{proof}
To begin, we know that $I \cdot I^{-1} = I$ so $I \in \G$ where
$I=[\eta_{ij}]_{n \times n}$ is the identity element of the group and satisfies
the condition $\eta_{ij} = 1$ every $i,j =1,\ldots,n$.

Now, observe that if $M \in \G$ then $M \cdot M^T = I$ and $M^T = M^{-1}$. Thus
$M^T(M^T)^T=I$ so $M^{-1} \in \G$.

Let $M$ and $N$ be arbitrary elements of $\G$. Observe that by the properties of $\G$:
$$NM(NM)^T=(NM)(N^T M^T) =N(MM^T) N^T=NIN^T=I.$$
%
$\G$ is closed and commutative under Hadamard product. Consequently, we see that $(\G,\cdot)$ is an abelian group.  
\end{proof}

\begin{df}
Let $\G$ be PC matrix Lie group and $M(t)$ be a path through $\G$. We say that $M(t)$ is \textit{smooth} if each entry in $M(t)$ is differentiable. The derivative of $M(t)$ at the point $t$ is denoted $M'(t)$ which is the matrix whose $ij^{th}$ element is the derivative of $ij^{th}$ element of $M(t)$.
\end{df}

\begin{cor}
The abelian group $\G$ is a PC matrix Lie group.
\end{cor}

\begin{proof}
We know that the Hadamard product "$\cdot$" and the operation $M\longrightarrow M^{-1}=M^T$ are differentiable for every PC matrix M.\\
Thus, $\G$ is a PC matrix Lie Group.
\end{proof}

We also know that the tangent space of any matrix Lie group at unity is a vector space. 

The tangent space of any matrix group $\G$ at unity $I$ will be denoted by $T_I(\G)=g$ where $I$ is the unit matrix of $\G$.

\begin{thm}
The tangent space of the PC matrix Lie Group $\G$ at unity $I$ consists of all $n \times n$ real matrices $X$ that satisfy $X + X^T = 0$.  
\end{thm}	
\begin{proof}
Recall that any matrix $A \in \G$ satisfies the condition $A\cdot A^T=I $.  Let us consider a smooth path $A(t)$ such that $A(0)=I$.
	
\noindent We know that:
	\begin{equation} \label{eq4-1}	
	A(t) \cdot A(t)^T = I 
	\end{equation}
for all parameters $t$. 
\noindent By differentiating the equation \ref{eq4-1}, we get
	\begin{equation} \label{4-2}
	 \frac{d}{dt}(A(t)) \cdot A(t)^T + A(t) \cdot \frac{d}{dt}(A(t)^T) = 0 
	\end{equation}
Considering that: $$ \frac{d}{dt}( A(t)^T) =  (\frac{d}{dt}(A(t)))^T $$ the equation \ref{4-2} 
can be rewritten as $$ A'(t) \cdot A(t)^T + A(t) \cdot
A'(t)^T = 0 $$
and at the point $t=0$, we obtain:
$$A'(0) + A'(0)^T = 0$$   
Thus, any tangent vector $X = A'(0)$ satisfies $X+X^T=0$.\\
\end{proof}
\begin{cor}
The Lie algebra of $\G$, denoted by $T_I(\G)$, is a Lie algebra of $\G$ and $T_I(\G)$ is the space of the skew-symmetric $n \times n$ matrices. Observe that:
$$dim(\G) = dim T_I(\G) = \frac{n \cdot (n-1)}{2}.$$
\end{cor}

\section{Exponential map and PC matrices}
\label{exp-trans}

The exponential map is a map from Lie algebra of a given Lie group to that group. In this Section, we will introduce the exponential transformation from $g$ (the tangent space to the identity element of PC matrix Lie group $\G$) to $\G$.\\

Let $\G$ be PC matrix Lie group and $g$ be the Lie algebra of $\G$. Then, the exponential map:
\begin{align*}
exp:g & \longrightarrow \G \\
A=[a_{ij}]_{n \times n} & \longrightarrow exp[A]=[e^{a_{ij}}]
\end{align*}

\noindent has the following properties: 
\begin{enumerate}
	\item $G_1=\{ \delta(t)=e^{tA} ~|~ t \in \mathbb{R}, A \in g \}$ is one parameter subgroup of $\G$.
	\item Let $A$ and $B$ be two elements of the Lie algebra $g$. Then, the following equality holds:
	$$e^{A+B} = e^A \cdot e^B$$
	\item Given any matrix $A \in g$, the tangent vector of the smooth path
	$\gamma (t)$ is equal to $A \cdot e^{tA}$, that is, $$\gamma '(t) =
	\frac{d}{dt}e^{tA} = A \cdot e^{tA} $$
	\item For any matrix  $A \in g$, 
          $$(e^A)^{-1} = e^{-A} = (e^A)^T = e^{A^T}$$
and
          $$ (e^A)(e^A)^T=e^A \cdot e^{A^T} = e^{A+A^T} = e^0 = 1.$$
	\item For any matrix  $A \in g$, $\gamma (t) = e^{tA}$ is a geodesic curve
of the pc matrix Lie group $\G$ passing through the point $\gamma (0) = 1$.
	
	\item For any matrix $A \in g$, we would like to stipulate that 
	$$det(e^A) = e^{Tr A}$$
	where $Tr(A)$ is the trace function of $A$.
	However, this cannot always be achieved and  a counterexample is presented in the Example~\ref{counterexample1}.  
	\end{enumerate} 

\begin{ex}
\label{counterexample1}
Let us consider the following matrix
$$ A=\begin{pmatrix}
		0 & -1 & 1 \\
		1 & 0  & 0 \\
		-1 & 0  & 0
	\end{pmatrix}
$$
		 
\noindent $A$ is the element of $g$, hence the 
exponential map of $A$ is:
$$ e^A=\begin{pmatrix}
1	& \frac{1}{e} & e \\
e	& 1 & 1 \\
\frac{1}{e} & 1 &  1
\end{pmatrix}
$$
		
\noindent The determinant of $e^A$ is: $$det(e^A) = e^2+e^{-2}-2$$ and the trace of $e^A$ is: $$Tr(A)= \Sigma^{3}_{i=1}a_{i,i} = 0$$ 

Consequently, $det(e^A)$ is not equal to $e^{Tr(A)}$ for the matrix $A$.  
\end{ex}

\section{Internal direct product of Lie group of $3 \times 3$ PC Matrices}
\label{dir-prod}

The aim of this section is to provide both geometric and algebraic perspectives on PC matrices. Our presentation is based on the techniques in \cite{A1995},\cite{YY2013}, and \cite{YY2016}. However, a modified approach is used in this Section. Let us recall that we consider only $3 \times 3$ PC matrices. Section~\ref{genSec} outlines generalization to $n>3$. Let us introduce the definition of the internal direct product.

\noindent We use the notation $I_n = \{1,\ldots n\}$.
\begin{df}
Let $\G$ be a group and let  $\{N_i|i\in I_n\}$ be a family of normal subgroups of $\G$. Then $\G=N_1N_2\dotsm N_n$ is called the {\em internal direct product} of $\{N_i | i \in I_n\}$
and $N_i\cap(N_1\dotsm N_{i-1}N_{i+1}\dotsm N_n)=\{e\}$ for all $i\in \{1,2,\ldots, n\}$ (see \cite{F1933}).
\end{df}

\begin{thm}
Let $\G$ be a group and $\{N_i | i \in I_n\}$ be a family of normal subgroups of $ \G $. Then $ \G $ is an internal direct
product of $\{N_i | i \in I_n\}$ if and only if for all $A \in \G$, $A$ can be uniquely expressed as $$A=a_1 \cdot a_2 \dotsm a_n$$ where $a_i \in N_i, i \in \{1,\ldots,n\}$ (see \cite{F1933})
\end{thm}


\begin{thm}
Let $\G$ be the internal direct product of a family of normal subgroups $\{N_i|i\in I_n\}$. 
Then $$\G \simeq N_1 \times N_2  \times \dotsm \times N_n$$.
		
The collection $M_n =(M_n,\cdot)$ of all consistent PC matrices
$ M $ is a multiplicative Lie subgroup of the Lie group $ \G $.
		 
Moreover, let us consider additive consistent matrices represented by the following set:
$$l = C_g = \{A = [a_{i,j}]_{n\times n}\in g | a_{i,k} + a_{k,j} = a_{i,j} \}
\text{ for every } (i,j,k) \in  T_n$$
\noindent where $T(n) = \{ (i,j,k) \in \{1,2,3, \dots ,n\}:i<j<k \}$ was 
considered in \cite{A1995}.

If $A, B \in C_g$ then $A+B \in C_g$ and $kA \in C_g$. Thus $C_g$ is a Lie
subalgebra of $(g,+)$. Let us observe that the following equality holds: 

\[g = C_g \oplus (C_g)^{\bot} (C_g)^{\bot} = h. \]

It follows that $C_{\G }$ is a Lie subgroup of $({\G},\cdot)$; therefore, the following equality also holds:
$$\G = C_{\G} (C_{\G})^{\perp}$$
\end{thm}  	
		
Considering the above results, we provide the new geometric and algebraic interpretation for $3 \times 3$ PC matrices.

\begin{prop}
Let
$$h=\Bigg\{\begin{pmatrix} 
0 & x & -x \\
-x & 0 & x \\
x & -x & 0
\end{pmatrix} : x \in \mathbb{R}\Bigg\}$$ 
and 
$$l=\Bigg\{\begin{pmatrix} 
0 & y & y+z \\
-y & 0 & z \\
-y-z & -z & 0
\end{pmatrix} : y,z \in \mathbb{R}\Bigg\}$$ 
be two sets of $3 \times 3$ PC matrices. Then the following holds:

\begin{enumerate}[(i)]
\item $h$ and $l$ are 1 and 2 dimensional Lie subalgebras of the Lie algebra of $g$,  
\item the vector space denoted by the Lie algebra of $h$ is the
orthogonal complement space of the vector space denoted by the Lie algebra of $l$.  
\end{enumerate}
\end{prop}
\begin{proof}
\begin{enumerate}[(i)]
	\item 
	For the proof, let us use  $sp$ for the linear span of a set of vectors.
	$$ h = sp \Bigg\{\begin{pmatrix} 
	0 & 1 & -1 \\
	-1 & 0 & 1 \\
	1 & -1 & 0
	\end{pmatrix}\Bigg\}$$ 
and
	$$ l = sp \Bigg\{\begin{pmatrix} 
	0 & 1 & 1 \\
	-1 & 0 & 0 \\
	-1 & 0 & 0
	\end{pmatrix},	\begin{pmatrix} 
	0 & 0 & 1 \\
	0 & 0 & 1 \\
	-1 & -1 & 0
	\end{pmatrix} \Bigg\}$$

Since the space $h$  produces one matrix, $dim(h) = 1$ and the space
$l$  produces two linearly independent vectors, $dim(l) = 2$.
Moreover, $h$ and $l$ are Lie subalgebras of $\G$.
	
	\item It is implied by the basic properties of Lie algebras. 
\end{enumerate}

We will show that every element of the Lie algebra $g$ can be written as the sum of one element of the Lie subalgebra $h$ and one element of Lie subalgebra $l$, that is, for all $A \in \G$, $A = A_h + A_l$, with $A_h \in h$ and $A_l \in l$.
$$
\begin{pmatrix} 
	0 & a & b \\
	-a & 0 & c \\
	-b & -c & 0
\end{pmatrix}
=
\begin{pmatrix} 
	0 & x & -x \\
	-x & 0 & x \\
	x & -x & 0
\end{pmatrix}
+
\begin{pmatrix} 
	0 & y & y+z \\
	-y & 0 & z \\
	-y-z & -z & 0
\end{pmatrix}
$$
where
$$x=\frac{1}{3}(a-b+c)$$
$$y=\frac{1}{3}(2a+b-c)$$
$$z=\frac{1}{3}(2c-a+b)$$
\end{proof}

\begin{prop}
\label{6-4}	
Let
$$H=\Bigg\{N|N=\begin{pmatrix} 
1 & k & \frac{1}{k} \\
\frac{1}{k}  & 1 & k \\
k & \frac{1}{k}  & 1
\end{pmatrix}, 
\mathit{where}\ k \in \mathbb{R^+}\Bigg\}$$ 
and 
$$L=\Bigg\{W|W=\begin{pmatrix} 
1 & y & yz \\
\frac{1}{y}  & 1 & z \\
\frac{1}{yz} & \frac{1}{z}  & 1
\end{pmatrix}, \mathit{where}\  y,z \in \mathbb{R^+}\Bigg\}$$ 

\noindent be two sets of $3 \times 3$ matrices. Then:

\begin{enumerate}[(i)]
	
   \item  $H$ and $L$ are normal subgroups of the Lie group of $\G$
   \item  The Lie group $\G$ is the internal direct product of the normal Lie
   subgroups  $H$ and $L$. In particular:
   $$\G \simeq H \times L.$$

\end{enumerate}
\end{prop}
\begin{proof}
\noindent We know that:

$$\G=\Bigg\{M|M=\begin{pmatrix} 
1 & m_{12} &  m_{13}\\
\frac{1}{m_{12}}  & 1 & m_{23} \\
\frac{1}{m_{13}} & \frac{1}{m_{23}}  & 1
\end{pmatrix} , m_{12},m_{13},m_{23} \in \mathbb{R^+}\Bigg\}$$ 

\ \\
Let us consider:\\
	$$A=\begin{pmatrix} 
1 & \xi & \eta \\
\frac{1}{\xi} & 1 & \gamma \\
\frac{1}{\eta}  & \frac{1}{\gamma}  & 1
\end{pmatrix} \in \G$$ 

\ \\
Using the logarithmic transformation we get:
$$\tilde{A}=\begin{pmatrix} 
0 & ln(\xi) & ln(\eta) \\
-ln(\xi) & 0 & ln(\gamma) \\
-ln(\eta) & -ln(\gamma) & 0
\end{pmatrix}$$

The Proposition~\ref{6-4} implies that: 

$$\tilde{A}=\begin{pmatrix} 
	0 & x & -x \\
	-x & 0 & x \\
	x & -x & 0
\end{pmatrix}
+
\begin{pmatrix} 
	0 & y & y+z \\
	-y & 0 & z \\
	-y-z & -z & 0
\end{pmatrix} = A_h+A_l$$

\noindent Moreover, we know that:
	$$x = \frac{ln(\xi)-ln(\eta)+ln(\gamma)}{3}= 
	      ln(\frac{\xi\gamma}{\eta})^{1/3}$$
	$$y = \frac{2ln(\xi)+ln(\eta)-ln(\gamma)}{3}= 
	      ln(\frac{\xi^2\eta}{\gamma})^{1/3}$$
	$$z = \frac{2ln(\gamma)-ln(\xi)+ln(\eta)}{3}= 
	      ln(\frac{\gamma^2\eta}{\xi})^{1/3}$$

\noindent Using the exponential transformation for $ \tilde{A} = A_h+A_l$	
we conclude that:
	$$ exp(\tilde{A}) = exp(A_h)\cdot exp(A_l)$$
which implies:
	$$A = A_H\cdot A_L$$
\end{proof}
	
Observe that the matrix $A_L$ is a consistent PC matrix. However, the following proposition needs to be considered for deciding on the classification of the matrix $A_H$.
    

\begin{prop}

\begin{enumerate}[(i)]
	\item Lie algebra of the PC matrix Lie group of consistent $L$ is $l$.
	\item Lie algebra of the PC matrix normal Lie subgroup $H$ is $h$.
\end{enumerate}
\end{prop}

 From the above proposition for additive consistent matrices, we conclude that $h$ is the orthogonal complement space of $l$. Hence every element of $h$
can be classified as ortho-additive consistent matrix and every element of $H$
is an ortho-consistent PC matrix.

\section{Generalization outline}
\label{genSec}
\begin{figure}[h]
	\includegraphics[width=1.2\linewidth]{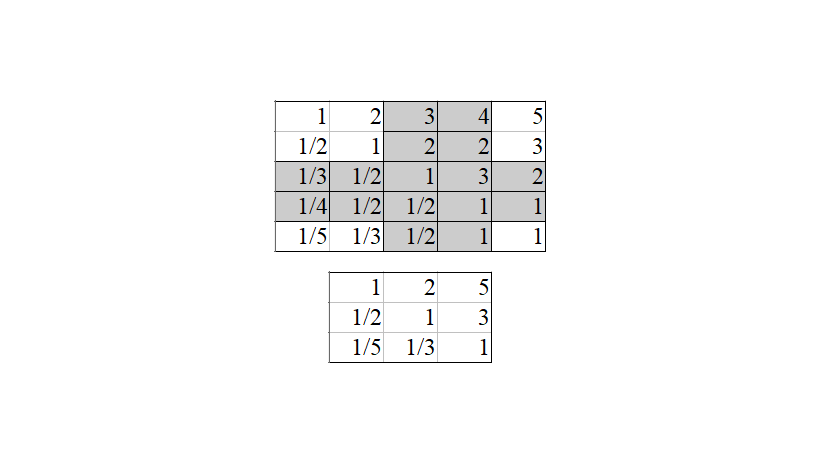}
	\caption[Generalization]{Generatization to $n>3$}
	\label{fig:m5tom3}
\end{figure}

Following \cite{KU2018}, we would like to point out that PC matrices of the size $n>3$ can be viewed as consisting of  $3 \times 3$ PC submatrices obtained by deleting rows and columns with the identical indices. That is, if we delete row
$j$, we need to also delete column $j$. Fig.~\ref{fig:m5tom3} demonstrates how to obtain one PC matrix $3 \times 3$ from PC matrix $5 \times 5$ by deleting
rows and columns numbered 3 and 4. 

Definition 3.3 from \cite{KU2018} states this as follows:
\begin{quotation}
	Assume $n<m$, $A$ and $B$ are square matrices of degree $n$ and $m$
	respectively. We call matrix $A$ a submatrix of $B$ ($A \subset B$) if there
	exist injection $$\sigma : \{1,\ldots,n\} \rightarrow \{1,\ldots,m\} $$ 
	such that for all 
	$n,m \in \{1,\ldots,n\}$ $a_{ij} = b_{\sigma(i)\sigma(j)}.$
\end{quotation}

The PC matrix reconstruction from its $3 \times 3$ submatrix approximation components is PC matrix $S$ of geometric means of all corresponding elements in these components. The need for geometric mean use comes from the occurrence of the same matrix element $(n-1)(n-2)(n-3)/6$ times in submatrices.

It is important to observe that PC matrix $S$ does not need to be consistent even though all submatrices of $S$
are consistent. However, the reconstruction process converges to a consistent PC matrix as partly proven by \cite{HK1996} and completed in \cite{KS2010}. The above reconstruction process will be analyzed in the planned follow-up paper.


\section*{Conclusions}

Using fundamental theorems from \cite{Levi1942} and \cite{Levi1943}, the collaborative research effort \cite{collab2020} demonstrated that group generalization for pairwise comparisons matrices is a  challenge. In particular, \cite{collab2020} provided the proof that elements of a multiplicative PC matrix must be selected from a torsion-free abelian group. 

Our paper demonstrates that the multiplicative PC matrices (not the elements of a PC matrix) generate a Lie group for the Hadamard product. Lie algebras of these Lie groups are identified here. It has been shown that Lie algebras form spaces of skew-symmetric matrices. It has also been proven that
the Lie group of PC matrices can be represented as an internal direct product using the direct summability property of vector spaces.

In conclusion, a relatively simple concept of pairwise comparisons turns out to be related to the theory of Lie groups and Lie algebras (what is commonly regarded as very sophisticated mathematics). For the first time, the decomposition of a PC matrix into an approximation component and orthogonal component (interpreted as the approximation error) was obtained. Without such decomposition, the pairwise comparisons method has remained incomplete for 722 years from its first scholarly presentation.

\section*{Acknowledgments}

We thank Tu\v{g}\c{c}e \c{C}alci for the verification of mathematical formulas and the terminology associated with them. The authors recognize the efforts of Tiffany Armstrong in proofreading this text and Lillian (Yingli) Song with the technical editing. We also acknowledge that algebraic terminology and basic concepts are based on \cite{Algebra}.


\begin{thebibliography}{99}
	
	\bibitem{A1995}
	Adukov, V., 
	On invertibility of matrix Wiener-Hopf operator on discrete linearly ordered abelian group, Integral Equations and Operator Theory, 23(4):373--386, 1995.
	
	\bibitem{YY2013} 
	Aksoyak, F.K.; Yayli, Y.,
	Homothetic Motions and Lie Groups in $\R_2^4$,
	Journal of Dynamical Systems and Geometric Theories,
	11(1--2): 23--38, 2013.
	
	\bibitem{AA2008}
	Arvanitoyeorgos, A.,
	Lie Transformation Groups and Geometry,
	In: Mladenov, I.M, (Ed),
	\textit{Ninth International Conference on Geometry, Integrability and Quantization},
	SOFTEX, 2008.
	
	\bibitem{BT2019}
	Bansal, S.; Tatu, A.,
	Affine Interpolation in a Lie Group Framework,
	ACM Transactions on Graphics, 38(4): 71, 2019. 
	
	\bibitem{BK1972}
	Beck, Re; Kolman, B.,
	Computer Approaches to Representations of Lie-Algebras,
	Journal of the ACM, 19(4):577-577, 1972. 
	
	\bibitem{YY2016}
	Bekar, M.; Yayli, Y.,
	Lie Algebra of Unit Tangent Bundle,
	Advances in Applied Clifford Algebras, 27(2): 965--975, 2017.
	
	\bibitem{BFK2011}
	Bozoki, S.; Fueloep, J.; Koczkodaj, W.W.,
	An LP-based inconsistency monitoring of pairwise comparison matrices,
	Mathematical and Computer Modelling, 54(1-2): 789-793, 2011.  
	
	\bibitem{CB2018}
	Cavallo, B.; Brunelli, M.,
	A general unified framework for interval pairwise comparison matrices, 
	International Journal of Approximate Reasoning, 93: 178--198, 2018.
	
	\bibitem{CavDap09}
	Cavallo, B.; D'Apuzzo L., A General unified framework for pairwise comparison matrices in multicriterial methods, International Journal of Intelligent Systems, 24(4): 377--398, 2009.
	
	\bibitem{RAND}
	Crawford, G.; Williams C., The Analysis of Subjective Judgment Matrices, 
	RAND Report R--2572--1--AF, 1985.
	
	\bibitem{D1999}
	Deng, HP,
	Multicriteria analysis with fuzzy pairwise comparison
	International Journal of Approximate Reasoning, 21(3): 215-231, 1999. 
	
	\bibitem{DAO2018}
	Demisse, G.G.; Aouada, D.; Ottersten, B..,
	IEEE Transactions on Pattern Analysis and Machine Intelligence,
	40(6): 1338--1351, 2018.
	
	\bibitem{HHPRD}
	Drton, M. ; H\"{a}gele, G.; Haneberg, D.; Pukelsheim, F. ; Reif  W.,
	The Augsburg Web Edition of Llull's Electoral Writings,
	\url{https://www.math.uni-augsburg.de/htdocs/emeriti/pukelsheim/llull/}
	
	\bibitem{EI2020}
	Ebisu, T.; Ichise, R.,
	Generalized Translation-Based Embedding of Knowledge Graph,
	IEEE Transactions on Knowledge and Data Engineering, 32(5): 941--951, 2020.
	
	\bibitem{F1933}
	Fitting, H, The theory of automorphism rings of Abelian groups and their analogue for non-commutative groups,	Mathematische Annalen, 107:514--542, 1933.
	
	\bibitem{HK1996}
	Holsztynski, W.; Koczkodaj, W.W., Convergence of inconsistency algorithms for the pairwise comparisons, Information Processing Letters, 59(4): 197--202, 1996.
	
	\bibitem{JK1996}
	Janicki, R; Koczkodaj, W.W.,
	A weak order approach to group ranking,
	Computers \& Mathematics with Applications, 32(2): 51-59, 1996.
	
	\bibitem{YY2014}
	Ismail G\"{o}k, I.; Okuyucu,  O.Z.; Ekmekci, N.;Yayli, Y.,
	On Mannheim partner curves in three dimensional Lie groups,
	Miskolc Mathematical Notes, 15(2): 467--479, 2014.
	
	\bibitem{WJ1991}
	Jefferys, W.H.; Berger, J.O., Ockham's Razor and Bayesian Statistics, American Scientist, 80(1): 64--72, 1991.
	
	\bibitem{KCD2009}
	Kobilarov, M.; Crane, K.; Desbrun, M.,
	Lie Group Integrators for Animation and Control of Vehicles,
	ACM Transactions on Graphics, 28(2): 16, 2009.
	
	\bibitem{K1996}
	Koczkodaj, W.W., 
	Testing the accuracy enhancement of pairwise comparisons by a Monte Carlo experiment, Journal of Statistical Planning and Inference, 69(1):21--31, 1996.
	
	\bibitem{K1998}
	Koczkodaj, W.W., 
	Statistically accurate evidence of improved error rate by pairwise comparisons,
	Perceptual and Motor Skills, 82(1):43--48, 1996.
	
	\bibitem{KSS2020}
	Koczkodaj, W.W.; Smarzewski, R.; Szybowski, J., On Orthogonal Projections on the Space of Consistent Pairwise Comparisons Matrices, Fundamenta Informaticae, 172(4): 379--397, 2020. 
	
	\bibitem{collab2020}
	Koczkodaj, W.W.; Liu, F; Marek, V.W. ; Mazurek, J; Mazurek, M; Mikhailov, L; Ozel, C; Pedrycz, W; Przelaskowski, A; Schumann, A; Smarzewski, R; Strzalka, 
	D; Szybowski, J; Yayli, Y,
	On the use of group theory to generalize elements of pairwise comparisons matrix: a cautionary note, International Journal of Approximate Reasoning, 124: 59--65, 2020.
	
	\bibitem{K2Y}
	Koczkodaj, W.W.; Mikhailov, L.; Redlarski, G.; Soltys, M.; Szybowski, J.; Tamazian, G.; Wajch, E.; Yuen, K.K.F., Important Facts and Observations about Pairwise Comparisons (the special issue edition), Fundamenta Informaticae, 
	144(3--4): 291--307, 2016.
	
	\bibitem{KO1997}
	Koczkodaj, W.W.; Orlowski, M,
	An orthogonal basis for computing a consistent approximation 
	to a pairwise comparisons matrix,
	Computers \& Mathematics With Applications, 34(10): 41--47, 1997.
	
	\bibitem{KS2010}
	Koczkodaj, W.W.; Szarek, S., 
	On distance-based inconsistency reduction algorithms for pairwise comparisons, 
	Logic Journal of the IGPL, 18(6): 859--869, 2010.
	
	\bibitem{KS2016}
	Koczkodaj, W.W.; Szybowski, J.,
	The Limit of Inconsistency Reduction in Pairwise Comparisons,
	International Journal of Applied Mathematics and Computer Science 26(3):721--729, 2016.
	
	\bibitem{KU2018}
	Koczkodaj, W.W.; Urban, R., 
	Axiomatization of inconsistency indicators for pairwise comparisons, International Journal of Approximate Reasoning, 94: 18--29, 2018.

\bibitem{KK2015}
Kulakowski, K., On the Properties of the Priority Deriving Procedure in the Pairwise Comparisons Method, Fundamenta Informaticae,139(4): 403-419, 2015.

	\bibitem{KMRS2019}
	Kulakowski, K.; Mazurek, J.; Ramík, J., Soltys, M.,
	When is the condition of order preservation met?
	European Journal of Operational Research, 277(1): 248--254, 2019.
	
	\bibitem{KLH2006}
	Kuo, M.-S.; Liang, G.-S.; Huang, W.-C.,		
	Extensions of the multicriteria analysis with pairwise comparison under a fuzzy environment
	International Journal of Approximate Reasoning, 43(3): 268--285, 2006.
	
	\bibitem{KF2020}
	Koudounas, A.; Fiori, S.,
	Gradient-based Learning Methods Extended to Smooth Manifolds Applied to Automated Clustering,
	Journal of Artificial Intelligence Research, 68: 777--815, 2020.
	
	\bibitem{LG2006}
	Lebanon, G.,
	Metric learning for text documents,
	IEEE Transactions on Pattern Analysis and Machine Intelligence, 28(4): 497--508, 2006. 
	
	\bibitem{Levi1942}
	Levi, F.W., Ordered groups, Proceedings of the Indian Academy of Sciences, Section A, 16: 256--263, 1942. 
	
	\bibitem{Levi1943}
	Levi, F. W., Contributions to the theory of ordered groups, Proceedings of the Indian Academy of Sciences, 17:199--201, 1943.
	
	\bibitem{LDX2019}
	Li, CC; Dong, YC; Xu, YJ; Chiclana, F.; Herrera-Viedma,; Herrera, F.,
	An overview on managing additive consistency of reciprocal preference relations for consistency-driven decision making and fusion: Taxonomy and future directions,
	Information Fusion, 52: 143--156, 2019.
	
	\bibitem{YY2010}
	Ozkaldi, S.; Yayli, Y.,
	Tensor Product Surfaces in $ \R^4 $ and Lie Groups,
	Bulletin of the Malaysian Mathematical Sciences Society, 33(1), 69--77, 2010.
	
	\bibitem{Ram16}
	Ram\'{\i}k, J., 
	Ranking alternatives by pairwise comparisons matrix with fuzzy element on Alo-group, in \textit{Intelligent decision technologies}, part II, book series: {\em Smart Innovation Systems and Technologies, Springer}, 57: 371--380, 2016.
	
	\bibitem{Algebra}
	Robinson, D.J.S., A Course in the Theory of Groups,
	Graduate Texts in Mathematics (GTM), vol. 80, Berlin, New York: Springer-Verlag, 1996.
	
	\bibitem{Tu2010}
	Tu, T.L., 
	An Introduction to Manifolds, 2nd ed., Springer, 2010.
	
	\bibitem{ZZH2019}
	Zhang, L.; Zheng, Q.-Z.; Huang, H., 
	Intrinsic Motion Stability Assessment for Video Stabilization,
	IEEE Transactions on Visualization and Computer Graphics,
	25(4): 1681--1692, 2019.	
	
	\bibitem{ZCX2019}
	Zhao, XL; Chen, QB; Xue, JF; Zhang, YM; Zhao, JJ,
	A Method for Calculating Network System Security Risk Based on a Lie Group, IEEE Access 7: 70610--70623, 2019.
	
\end{thebibliography}
\end{document}